\documentclass[12pt]
{article}
\usepackage{amssymb,amsmath, mathtools}

\usepackage[latin1]{inputenc}
\usepackage[T1]{fontenc}
\usepackage{play}
\usepackage{mathtools}
\usepackage{graphicx}
\usepackage{color}
\usepackage{booktabs,multirow}
\usepackage{caption}
\usepackage{subcaption}

\usepackage[table,xcdraw]{xcolor}
\DeclareGraphicsExtensions{.pdf,.png,.jpg}

\DeclareMathOperator{\sgn}{sgn}

\graphicspath{ {images/} }

\usepackage{lipsum}

\let\OLDthebibliography\thebibliography
\renewcommand\thebibliography[1]{
  \OLDthebibliography{#1}
  \setlength{\parskip}{0pt}
  \setlength{\itemsep}{0pt plus 0.3ex}
}

\usepackage{amsmath}
\usepackage{amsthm}
\usepackage{amsfonts}
\usepackage{bbm}
\usepackage{amssymb}
\usepackage{graphicx}
\usepackage{color}

\usepackage{wasysym, marvosym}
\usepackage{lmodern}

\usepackage{fancybox}
\newlength{\querylen}
\newcommand{\prob}{\mathbb{P}}

\newcommand{\me}{\mathbb{E}}

\newcommand{\profile}{\boldsymbol p}

\newcommand{\prior}{\boldsymbol \pi}

\newtheorem{thm}{Theorem}
\newtheorem{lemma}[thm]{Lemma}

\newtheorem{assertion}[thm]{Proposition}
\theoremstyle{definition}

\theoremstyle{remark}

\begin{document}
\title{The  Last-Success Stopping Problem with Random Observation Times
} 
\author{Alexander Gnedin ~~~and~~~Zakaria Derbazi\\ {\it\small Queen Mary, University of London}}

\maketitle
\noindent

\begin{abstract}
\noindent
Suppose $N$ independent Bernoulli trials  are observed sequentially 
at  random times of a mixed binomial process. 
The task is to maximise, by using a nonanticipating stopping strategy, the probability of stopping 
 at   the last success. 
We focus on  the version of the problem where the $k^\text{th}$ trial is a success with probability $p_k=\theta/(\theta+k-1)$ and the prior distribution of $N$ is negative binomial  with shape parameter $\nu$.
Exploring  properties of the Gaussian hypergeometric function, we find that 
the myopic stopping strategy
is optimal if and only if  $\nu\geq\theta$. 
We derive formulas to assess the winning probability and discuss limit forms of the problem for large $N$.

\end{abstract}

\section{Introduction}
Given   $N$ sequentially observed independent Bernoulli trials with  success probability  $p_k$ for  the $k^\text{th}$ trial,
the last-success problem  asks one to find a stopping strategy  maximising  the probability of stopping  
 at the last success. The  instance  $p_k=1/k$ corresponds  to the famous best-choice problem, in 
which the observer compares  rankable items arriving in random order   with the goal to stop at the best \cite{Handbook}.
If $N=n$, is a fixed  number known to the observer, then the optimal strategy has a particular simple form:  there exists a
  critical index $k_n$, such that one  should always skip  the initial
$k_n$ trials then stop at the first consequent success
 \cite{Odds, HillKrengel,  Ribas}. 
In this paper we will consider an instance of   the following more complex setting,
which generalises some previously studied best-choice problems with 
random  $N$ and random observation times
 {\cite{Browne,  BrussCZ, Bruss88, BrussSam1, BrussSam2, CZ,  Gaver, AnoK,  Stewart, TamakiWang}}.

 Let the number of trials $N$ be governed by the distribution  $\prior=(\pi_0,\pi_1,\cdots)$, called  in the sequel {\rm prior}, 
and let $\profile=(p_1,p_2,\cdots)$ be a profile of success probabilities. Suppose
the trials are paced randomly in time according to a mixed binomial point process with some continuous base distribution.
That is to say, conditionally on $N=n$  the epochs of the trials  are allocated like $n$  order statistics from the base distribution.
The trial at the $k^\text{th}$ epoch
is a success with probability $p_k$, independently of the outcomes and epochs of all other trials.
The base distribution, $\prior$ and $\profile$ are assumed to be known to the observer, who learns the times and outcomes of the trials online, 
aiming to stop, by means of  a nonanticipating stopping strategy, at the moment of the last success.
The form of the base distribution is not important, and  without losing generality in the sequel we shall assume it  to be  $[0,1]$-uniform. 

In the best-choice problem with random observation times, a distinguished role is played by the so-called $1/e$-strategy,  which prescribes to skip all trials before time $1/e$ then stop at the first 
success.   Bruss \cite{Bruss84} observed 
that  for $N\geq 1$ the $1/e$-strategy  guarantees winning probability $1/e$ and  that this bound is  asymptotically sharp for $N$ large.
However, for any fixed prior $\prior$ the $1/e$-strategy is not optimal \cite{ G1/e}.

The analogy with many other problems of optimal stopping \cite{Ferguson} suggests
in search for the optimum to  examine the {myopic} strategy.   This strategy, denoted here $\tau^*$, prescribes to stop as soon as  the probability to win  
 by stopping at the current trial  is at least  as great as the  probability to win by stopping at the next success (if available).
Concretely,  if  the  $k^\text{th}$ trial observed at time $t$ is a success, this stopping condition holds whenever $t\geq a_k$ for certain cutoff $a_k$. 
Optimality of the  myopic strategy is ensured whenever  the sequence of cutoffs is nonincreasing. Here is a summary of what has been known to date about the optimality 
of $\tau^*$ in the context of the best-choice problem:

\begin{itemize}
\item[--]  Cowan  and Zabczyk \cite{CZ}  proved optimality of  $\tau^*$ for the Poisson prior (also see \cite{BG,  BrussSam1, Ciesielski} for asymptotics 
of the cutoffs). 
\item[--] Tamaki and Wang \cite{TamakiWang} did this for the discrete uniform distribution.
\item[--] Bruss \cite{BrussCZ} found that  for the geometric prior, $\tau^*$ is optimal and has a single cutoff (i.e. the $a_k$'s are all equal).
This is a paradigmatic case, characterised by the property that the epochs of success trials comprise a Poisson process
\cite{BrowneBunge, BrussRogers}. 
\item[--] The problem for the  negative binomial prior  with integer shape parameter $\nu>1$ 
was first stated  in \cite{BrussCZ} in terms of a mixed Poisson representation of the pacing process.
 The optimality of $\tau^*$ was claimed
 in an unpublished paper  by Browne \cite{Browne}, but  the crucial monotonicity property  of the cutoffs  was left  without proof.
The latter was settled  independently in an analytic tour de force 
by  Kurushima and Ano \cite{AnoK},  who also used the mixed Poisson representation.
\item[--] For the logarithmic series prior we have shown recently that  the cutoffs are increasing,  and that $\tau^*$ is not optimal and  can be improved by some explicitly constructed strategies \cite{GD}. 
In different terms, this instance of the problem had been introduced in  Bruss and Yor \cite{BrussYor} as a way to model the observer
with  complete ignorance about $N$.
Bruss and Rogers  \cite{BrussRogers21} identified the structure of the observation process and disproved the conjecture 
about optimality of the $1/e$-strategy.
\end{itemize}

 A central theme in this paper is the last-success problem in 
 the setting where  $\prior$ is a negative binomial  distribution  with arbitrary shape parameter $\nu>0$, and where the  profile of success probabilities is
\begin{equation}\label{profile}
p_k=\frac{\theta}{(\theta+k-1)}\,, ~~~k\geq 1,
\end{equation}
with another parameter     $\theta>0$. The case $\theta=1$ corresponds to the best-choice problem studied in \cite{BrussCZ} for  $\nu=1$ and  in  \cite{Browne, AnoK} for integer $\nu>1$.
The  more general profile (\ref{profile}) appears in connection with a
 nonstationary model of the  theory of records \cite{Nevzorov, Pfeifer}, so from this perspective
 the  paper falls in the area of random record models; see \cite{BrowneBunge, Gaver, Orsingher, OB} and Section \ref{Nevzorov} of this paper.
For integer $\theta$, the last-success problem can be interpreted as  the best-choice problem with batch observations \cite{Govindarajulu, Hsiau}.
Our main motivation, however, stems from  connections to the 
broad  circle of questions around the famous Ewens sampling formula \cite{ABT, Crane}. The latter  
suggests a  species  sampling metaphor for the following problem:   given a time limit, how to  stop exploring a population comprised of distinct unknown species
at the very moment when the last new species is discovered?

We  will show that the myopic strategy is optimal for $\theta\geq\nu$ and not optimal for $0\leq \nu<\theta$, where the edge case $\nu=0$ corresponds to the logarithmic series prior.
The direction of monotonicity of the cutoffs sequence changes at $\theta=\nu$.
The approach taken relies on  properties of the Gaussian hypergeometric function (summarised in Section \ref{HGF})   that appear to be  of independent interest.
In the watershed case $\theta=\nu$, the success epochs 
 comprise a Poisson process and all cutoffs are equal, similarly to 
the instance $\theta=\nu=1$ from \cite{BrussCZ}. 
A counterpart of the best-choice problem treated in \cite{Browne, AnoK} is the special case $\nu-\theta\in{\mathbb N}$, where 
the relevant hypergeometric series are polynomials.


To assess the strategies with nonincreasing cutoffs, we derive a formula for the winning probability, very much in spirit of 
the formulas known for the full information best-choice problem \cite{GM, Poros1}. 
We also discuss some asymptotic theory for large $N$.
For the  profile (\ref{profile}), the limit winning probability is still $e^{-1}$,  achieved by a strategy with single cutoff 
 $e^{-1/\theta}$.
A more delicate  limit    model  
exploits the idea of  an infinite prior, as appeared    under different guises in the work on the best-choice problem 
\cite{BrussCZ,   Bruss88, BrussRogers, BrussRogers21, BrussSam1, GD,  Sakaguchi, Stewart}.


\section{The problem with fixed $N=n$}\label{S2}

In this section we  sketch  the last-success problem  where the number of trials $n$  is fixed and known to the observer \cite{Odds, GD, HillKrengel, Ribas}.
Consider first the general profile  $\profile$ of success probabilities, with $p_1=1$ and  $0<p_k<1$ for $k>1$.
The number of successes among trials $k+1,\cdots,n$  has  the Poisson-binomial distribution with probability generating function
$$z\mapsto \prod_{j=k+1}^n (1-p_j+zp_j).$$
In particular, the  probability of  zero successes is
\begin{equation}\label{s0}
s_0(k+1, n):=\prod_{j=k+1}^n (1-p_j),~~~~k\geq0,
\end{equation}
and the probability of one success is
\begin{equation}\label{s0s1}
s_1(k+1, n):=s_0(k+1,n) \sum_{j=k+1}^n \frac{p_j}{1-p_j}, ~~~~~k\geq 1,
\end{equation}
where
$n\geq k$ and $s_1(k+1,k)=0$.
Since the sum decreases in $k$ 
there exists a threshold $k_n$ such that
the inequality
$$s_0(k+1,n)\geq  s_1(k+1,n)$$
holds for $k\geq k_n$. 
It follows as in the classic best-choice problem \cite{Dynkin} that the optimal strategy of the observer knowing $n$ is to stop  at the first success among trials $k_n+1,\cdots,n$ (or to not stop at all if none
of these trials is a success).
The function $s_1(k+1,n)$ is unimodal in $k$, with maximum achieved at the threshold value $k_n$, see
\cite{Odds, GD, HillKrengel, Pfeifer, Ribas}  for
more details  including asymptotics and bounds.

It is also useful to inspect the dependence on $n$.
As $n$ grows, both $k_n$ and the maximum value $s_1(k_n+1,n)$ are nonincreasing, provided $p_{k_n+1}\geq p_{n+1}$ (this condition certainly holds if  $p_k\downarrow$), see \cite{GD}. Moreover, we have:

\begin{lemma}\label{LL1} The function $s_1(k+1,k+j)$  is unimodal in the variable $j$.
\end{lemma}
\begin{proof} Using (\ref{s0s1}) and manipulating the sums gives
\begin{eqnarray*}
s_1(k+1,k+j+1)-s_1(k+1,k+j)=s_0(k+1, k+j) p_{k+j+1}\left(  1-\sum_{i=k+1}^{k+j} \frac{p_i}{1-p_i} \right).
\end{eqnarray*}
Since the expression in brackets decreases in $j$, the left-hand side has at most one change of sign from $+$ to $-$.
\end{proof}

Now, we specialise the success profile to (\ref{profile}).
The distribution of the number of successes among trials $k+1,\cdots, n$  has a nice  probability generating function
\begin{equation}\label{PGF}
z\mapsto \prod_{i=k+1}^n  \left(1-\frac{\theta}{\theta+i-1}+ \frac{\theta z}{\theta+i-1}   \right)   = \frac{(k+\theta z)_{n-k}}{(k+\theta)_{n-k}},
\end{equation}
where $(x)_n$ denotes the Pochhammer factorial.
The instance
$k=0$  is sometimes called the Karamata-Stirling law \cite{Bingham}.
From (\ref{PGF}),
for $k\geq 1$
the probability of zero and one successes becomes, respectively,
$$s_0(k+1,n)=  \frac{(k)_{n-k}}{(k+\theta)_{n-k}} ~~~{\rm and}~~~s_1(k+1,n)=  s_0(k+1,n) \sum_{i=k+1}^n \frac{\theta}{i-1}.$$
The threshold $k_n$  is found from the double inequality
\begin{equation}\label{d-ineq}
\frac{1}{k_n-1} +\cdots+\frac{1}{n-1} \leq \frac{1}{\theta}   < \frac{1}{k_n} +\cdots+\frac{1}{n-1}\,,
\end{equation}
in full analogy with the classic case $\theta=1$    \cite{Dynkin}.

As $n\to\infty$, the logarithmic approximation to the harmonic series yields via  (\ref{d-ineq}) the  limits 
\begin{equation}\label{fixed-n}
 \frac{k_n}{n}\to e^{-1/\theta}~~~{\rm and~~~} s_1(k_n+1,n)\downarrow e^{-1}.
\end{equation}
Thus  the limiting winning probability is $e^{-1}$ regardless of $\theta$, but the proportion of trials skipped before being willing to stop on a success depends on the parameter.

\section{Myopic and other strategies in the general setting}

Consider  the general prior $\prior$ and profile $\profile$.
Let $N, U_1, U_2,\cdots$ be independent random variables,  with $U_n$'s being   uniform-$[0,1]$ and $\prob(N=n)=\pi_n$.
Define a mixed binomial process 
$$
N_t=\sum_{n=1}^{N}  1_{\{U_n\leq t\}}, ~~~t\in [0,1],
$$
with chronologically ordered epochs $T_1<T_2<\cdots<T_{N}$. That is, given $N=n$ (for $n\geq1$)      the epochs  $T_1,\cdots, T_n$ coincide with the order statistics of $U_1,\cdots,U_n$.
With each $T_k$ we associate a Bernoulli trial resulting in success with probability $p_k$, independently of anything else.
We let $T_n$ undefined on the event $\{N<n\}$.

We speak of  {state} $(t,k)\in [0,1]\times\{0,1,\cdots\}$   to designate the event $\{N_t=k\}$.
Similarly,  under state $(t,k)^\circ$  (for $k\geq1$)  
we mean that the time of trial with index $k$ is  $T_k=t$   and that the trial is a success. 
Being  a subsequence of $(T_1,1),\cdots,(T_N,N)$, the bivariate sequence of success trials  increases in both components.
The generic success  trial can also be denoted as $(t, N_t)^\circ$, meaning that $t$ is some unspecified $T_k$ and the $k^\text{th}$ observed  trial  is a success.

Using the term `state' is justified by  the fact  that the sequence of success trials is Markovian, because 
$(t, N_t)$ is  all what we need to know  to make  inference about the number  of  trials after time $t$ and their outcomes.
Formally, let $({\cal F}_t, ~t\in[0,1])$ be the natural filtration, where the $\sigma$-algebra ${\cal F}_t$ embraces all the information about epochs and outcomes of the trials on $[0,t]$.
Then by  independence of the trials and the
Markov property of  mixed binomial processes \cite{Kallenberg},
 the distribution of the point process of epochs and outcomes of the trials on $(t,1]$ both  depend on ${\cal F}_t$ only through $N_t$.

The posterior distribution of the number of trials  factors as 
\begin{eqnarray}\label{posterior}
\prob(N=k+j\,|\,N_t=k)=f_k(t) {k+j\choose k} \pi_{k+j} (1-t)^j,~~~j\geq 0,
\end{eqnarray}
where $f_k(t)$ is a normalisation function.
By (\ref{posterior}) the conditional probability that there are no successes following state  $(t,k)$ can be written as
\begin{equation}\label{S0}
{\cal S}_0(t,k):=f_k(t) \sum_{j=0}^\infty \pi_{k+j} {k+j\choose k}  (1-t)^j s_0(k+1, k+j),
\end{equation}
and the probability that there is exactly one success as
\begin{equation}\label{S1}
{\cal S}_1(t,k):=f_k(t) \sum_{j=1}^\infty \pi_{k+j} {k+j\choose k}  (1-t)^j s_1(k+1, k+j),
\end{equation}
where $s_0$ and $s_1$ are defined by (\ref{s0}) and (\ref{s0s1}), respectively.
The {\it adapted} reward 
from stopping at $(t,k)^\circ$ is the conditional winning probability  equal to ${\cal S}_0(t,k)$.
Likewise, the adapted reward from stopping at the next success is  the probability   ${\cal S}_1(t,k)$.

\begin{lemma}\label{Lemma1} For fixed $k\geq1$, suppose $\pi_k>0$. Then 
\begin{itemize}
\item[\rm(i)]  ${\cal S}_0(1,k)=1, ~{\cal S}_1(1,k)=0$,
\item[\rm(ii)]   ${\cal S}_0(t,k)$ is nondecreasing in $t$, and is strictly increasing   if   $~\sum_{j=k+1}^\infty \pi_j>0$,
\item[\rm(iii)]   
there exists $a_k\in[0,1)$ such that 
$\sgn ({\cal S}_0(t,k)-{\cal S}_1(1,k)) =\sgn (t-a_k)$ for  $t\in(0,1]$.
\end{itemize}
 \end{lemma}
\begin{proof} Part (i) follows by checking that $\prob(N=k\,|\,N_t=k)\to 1$ as $t\to1$.

For (ii) we first show that   $N$ conditioned on $N_t=k$ stochastically increases with $k$.
To that end, consider two dependent  copies $(M_t), (N_t)$ of the same Markov counting process.
Start $(M_t)$ at $(t_0, k+1)$ and $(N_t)$ at $(t_0, k)$ and let them run independently until  random time $\rho$, defined to be the time when 
the processes hit the same value, so $M_\rho=N_\rho$, or $\rho=1$ if   this does not occur. In the event $\rho<1$  let the processes get coupled from time $\rho$ on.
Since both processes have only unit upward jumps, $M_1>N_1$ in the event $\rho=1$, while $M_1=N_1$ in the event $\rho<1$, so in any case $M_1\geq N_1$.
For $N=N_1$ this obviously gives 
$$\prob(N\geq n\,|\,N_t=k)\leq \prob(N\geq n\,|\,N_t=k+1),$$
which is the desired stochastic monotonicity.
Next, since $(N_t)$ is nondecreasing and Markovian, for $s<t$ the distribution of $N-N_s$ given $N_s=k$ dominates stochastically the distribution of $N-N_t$ given $N_t=k$.
Let $J_1$ have the first of these two distributions and $J_2$ the second. 
As $J_1$ is stochastically larger than $J_2$, the number of successes among trials $k+1,\cdots,k+J_1$  is stochastically larger than the number of successes among trials
$k+1,\cdots,k+J_2$. In particular, the probability of zero successes is smaller for trials $k+1,\cdots,k+J_1$ than for trials $k+1,\cdots,k+J_2$, which is
  ${\cal S}_0(s,k)\leq {\cal S}_0(t,k)$. The inequality is strict if the distributions of $J_1$ and $J_2$ are distinct, which is ensured by the condition $\prob(N>k)=\sum_{j=k+1}^\infty\pi_j>0$.

By (\ref{s0}) and (\ref{s0s1}),  $s_0(k+1,k+j)-s_1(k+1, k+j)$ has at most one change of sign as $j$ increases. 
Applying Descartes' rule of signs    for power series \cite{Descartes}, ${\cal S}_0(t,k)-{\cal S}_1(t,k)$ has at most one root on $[0,1)$.
In the view of (i),
 either the difference is positive everywhere on $(0,1]$ or the sign switches from $-$ to $+$ at some $a_k\in (0,1)$.
\end{proof}

We define a {\it stopping strategy} to be  a random variable $\tau$ which takes  values in the random set of times $\{T_1,\cdots, T_N, 1\}$ and is adapted in the sense that  $\{\tau\leq t\}\in {\cal F}_t$ for  $t\in[0,1]$.
The event $\{\tau=1\}$ is interpreted as not stopping at all. 
We shall restrict consideration to the strategies that only stop at successes $(T_k, k)^\circ$.
Under stopping in state
 $(t,k)^\circ$ we mean $\tau=T_k$ and that $T_k=t$ is a success epoch.
For $k\geq 1$ stopping at  $(t,k)^\circ$  is regarded as a win if  $T_k=t$ is the last success epoch.

With every `stopping set'  $B\subset [0,1)\times{\mathbb N}$ we associate a {\it Markovian} stopping strategy 
which stops at the first success trial 
$(t,k)^\circ$  with  $(t,k)$ falling in  $B$.
The optimal strategy exists for arbitrary $\prior$ and $\profile$,
and is Markovian by a theorem on exclusion of randomised stopping strategies (see \cite{CRS}, Theorem 5.3).

The {\it myopic} strategy $\tau^*$ is defined as the strategy stopping at the first state $(t,k)^\circ$ satisfying ${\cal S}_0(t,k)\geq {\cal S}_1(t,k)$.
This strategy is Markovian, with the stopping set being
\begin{equation}\label{A}
A:=\bigcup\limits_{k=1}^\infty ([a_k,1)\times \{k\}),
\end{equation}
where  the cutoff  $a_k$ is defined in Lemma \ref{Lemma1};
 this is the earliest time when $\tau^*$ can accept  success trial with index $k$. 
If  $(t,k)\notin A$ then it is not optimal to stop in the state  $(t,k)^\circ$; but the converse is not true in general.

The so-called {\it monotone case}   of  the optimal stopping theory holds if  $A$ is 
`closed',  meaning that after entering $A$, the sequence of successes  $(t,N_t)^\circ$  does not exit the set \cite{CRS, Ferguson}.
This condition is necessary and sufficient for the optimality of  $\tau^*$.
Since the sequence of successes is increasing in both components,
this property of $A$ is equivalent to the condition
$a_1\geq a_2\geq\cdots$, which is characteristic for the optimality  of $\tau^*$.



\section{Power series priors}
A distribution of $N$ can be treated within a family of power series priors 
\begin{equation}\label{p-series}
\pi_n= c(q) w_n q^n, ~~~n\geq 1,
\end{equation}
which share the same shape  $(w_n)$ and have scale parameter  $q>0$ varying within the radius of convergence of  
 $\sum_n w_n q^n$.
A rationale for such extension is  that the last-success problems are consistent for  different $q$, 
since   the posterior distribution of $N$ in state $(t,k)$ depends on $t$ through the variable
\begin{equation}\label{cd}
x=q(1-t).
\end{equation}
Explicitly,
 (\ref{posterior}) becomes
\begin{eqnarray*}
\prob(N=k+j\,|\,N_t=k)&=&[f_k(t) c(q)q^k] {k+j\choose k} w_{k+j} [q(1-t)]^j\\
&=&
\tilde{f}_k(x)    {k+j\choose k} w_{k+j} x^j,
\end{eqnarray*}
where  $\tilde{f}_k(x)$ is a normalisation function.
Casting (\ref{S0}) and (\ref{S1}) as power series in $x$,
the cutoffs can be defined in terms of 
a sole sequence of critical roots $(\alpha_k)$ not depending on $q$  via the formula
\begin{equation}\label{a-alpha}
a_k=\left(1-\frac{\alpha_k}{q} \right)_+.
\end{equation}
The myopic strategy is optimal for all $q$ if and only if  $\alpha_k\uparrow$, that is the roots are nondecreasing with $k$.

\section{The problem under the negative binomial prior}

We proceed  with the last-success problem for the success probabilities  (\ref{profile}) and $N$ following the negative binomial prior ${\rm NB}(\nu,q)$  with shape parameter $\nu>0$ and scale parameter $0<q<1$. 
\begin{equation}\label{NB}
\pi_n=     \frac{(\nu)_n}{n!}  (1-q)^\nu q^n, ~n\geq 0,
\end{equation}

The distribution increases stochastically in both $q$ and $\nu$, as the formula for the mean
$
\me N= {q\nu}/{(1-q)}
$
suggests.

The pacing process has a nice structure:
\begin{assertion}\label{Structure} Under the negative binomial prior ${\rm NB}(\nu,q)$ it holds that 
\begin{itemize}

\item[\rm(i)] the distribution of $N_t$ is ${\rm NB}(\nu,qt/(1-q+qt))$,
\item[\rm(ii)] the posterior of $N-N_t$ given $N_t=k$ is ${\rm NB}(\nu+k, q(1-t))$,

\item[\rm(iii)] $(N_t,~t\in[0,1])$ is a P{\'o}lya-Lundberg  birth process  with jump rate 
\begin{equation}\label{JR}
\frac{k+\nu}{t+q^{-1}-1}
\end{equation}
in state $(t,k)$,

\item[\rm(iv)]  $(N_t,~t\in[0,1])$ is a mixed Poisson process with random rate distributed according to
 ${\rm Gamma}(\nu,q^{-1}-1)$. 
\end{itemize}

\end{assertion}

\begin{proof}

Part (ii)   follows from (\ref{posterior}), (\ref{NB})  and the identity
$${k+j\choose k} \frac{(\nu)_{k+j}}{(k+j)!}=\frac{(\nu)_k}{k!}                         \frac{ (\nu+k)_j}{j!}.$$
The posterior mean of the number of trials on $[t,1]$ is thus
$$\frac{q(1-t)(\nu+k)}{1-q(1-t)},$$
and the formula for the rate in part (ii) follows since the epochs of the trials are spread by the $[t,1]$-uniform distribution.
Part (iv) follows from the mixture representation of ${\rm NB}(\nu,q)$.  Other assertions are straightforward.
\end{proof}
\noindent
See \cite{Kotz} for summary of properties and applications of the P{\'o}lya-Lundberg process.


Let  $F(a,b,c,x)$ be the Gaussian hypergeometric function (commonly denoted $_2F_1(a,b;c; x)$).
The properties we need to analyse $\tau^*$ are collected in Section \ref{HGF}.
It will be convenient to use interchangeably the variables 
$x=q(1-t)$ and $t=1-x/q.$

From  (\ref{PGF}) the probability generating function of the number of successes following state $(t,k)$ becomes
\begin{eqnarray}\nonumber
z \mapsto(1-x)^{k+\nu}\,\sum_{j=0}^\infty \frac{(k+\nu)_j}{j!}\frac{(k+\theta z)_j}{(k+\theta)_j}\,  x^j&=& (1-x)^{k+\nu} F(k+\theta z, k+\nu, k+\theta; x)\\
\nonumber
&=&
(1-x)^{\theta-\theta z}F(\theta-\theta z,\theta-\nu,\theta+k,x), \\
\label{PGF1}
\end{eqnarray}
where the second representation
follows from 
the identity (\ref{change-par}). 
Thus the probability of
 zero, respectively one successes is given by
\begin{eqnarray*}
{\cal S}_0(1-x/q,k)&=&(1-x)^{k+\nu}F(k, k+\nu, k+\theta; x)\\
&=&
(1-x)^\theta F(\theta,\theta-\nu, k+\theta,x).\\
		{\cal S}_1(1-x/q,k)&=&\theta (1-x)^{k+\nu}D_a F(k, k+\nu, k+\theta; x)\\
&=& 
-\theta (1-x)^\theta
\left\{ \log (1-x) F(\theta,\theta-\nu,\theta+k,x) 
+D_a  F(\theta,\theta-\nu,\theta+k,x) \right\}
\end{eqnarray*}
where $D_a$ stands for the derivative in the first parameter of $F(a,b,c,x)$.  

Complementing the general Lemma \ref{Lemma1} we have specifically for the setting in focus:

\begin{lemma}$   $

\begin{itemize}
\item[\rm (i)] ${\cal S}_0(t,k)$ is increasing in $k$ for  $\nu>\theta$, and decreasing in $k$ for $\nu<\theta$,
\item[\rm (ii)] ${\cal S}_1(t,k)$  is unimodal in $t$.
\end{itemize}
\end{lemma}

\begin{proof} Part (i) follows from Lemma \ref{L0-HH} applied to $F(\theta,\theta-\nu,k+\theta)$.
$${\cal S}_1(1-x/q,k)=(1-x)^{k+\nu} \sum_{j=0}^\infty   \frac{(\nu+k)_j}{j!} \, x^j  \,      s_1(k+1, k+j)            $$ 
Part (ii) follows from the unimodality of $s_1(k+1, k+j)$ in $j$ (Lemma \ref{LL1})  and the shape-preserving property in Lemma \ref{Baskakov}.
\end{proof}

We introduce the notation $W_0(x, k) = {\cal S}_0(1-x/q,k)$ and $W_1(x, k) = {\cal S}_1(1-x/q,k)$ 
 to stress that these two probabilities do not really depend on $q$.

For $k\geq 1$, equation
$$W_0(x, k)= W_1(x, k)$$
becomes
\begin{equation}\label{Eq-alpha}
D_a \log F(\theta,\theta-\nu,\theta+k,x)=-\frac{1}{\theta} -\log(1-x).
\end{equation}
Define  $\alpha_k$ to be the (unique)  root of the equation on the unit interval.

Denote 
$$\alpha^*:= 1-e^{-1/\theta}.$$

\begin{figure}[h]
	\centering
	\begin{subfigure}[b]{0.5\textwidth}
		\centering
		\includegraphics[width=\textwidth]{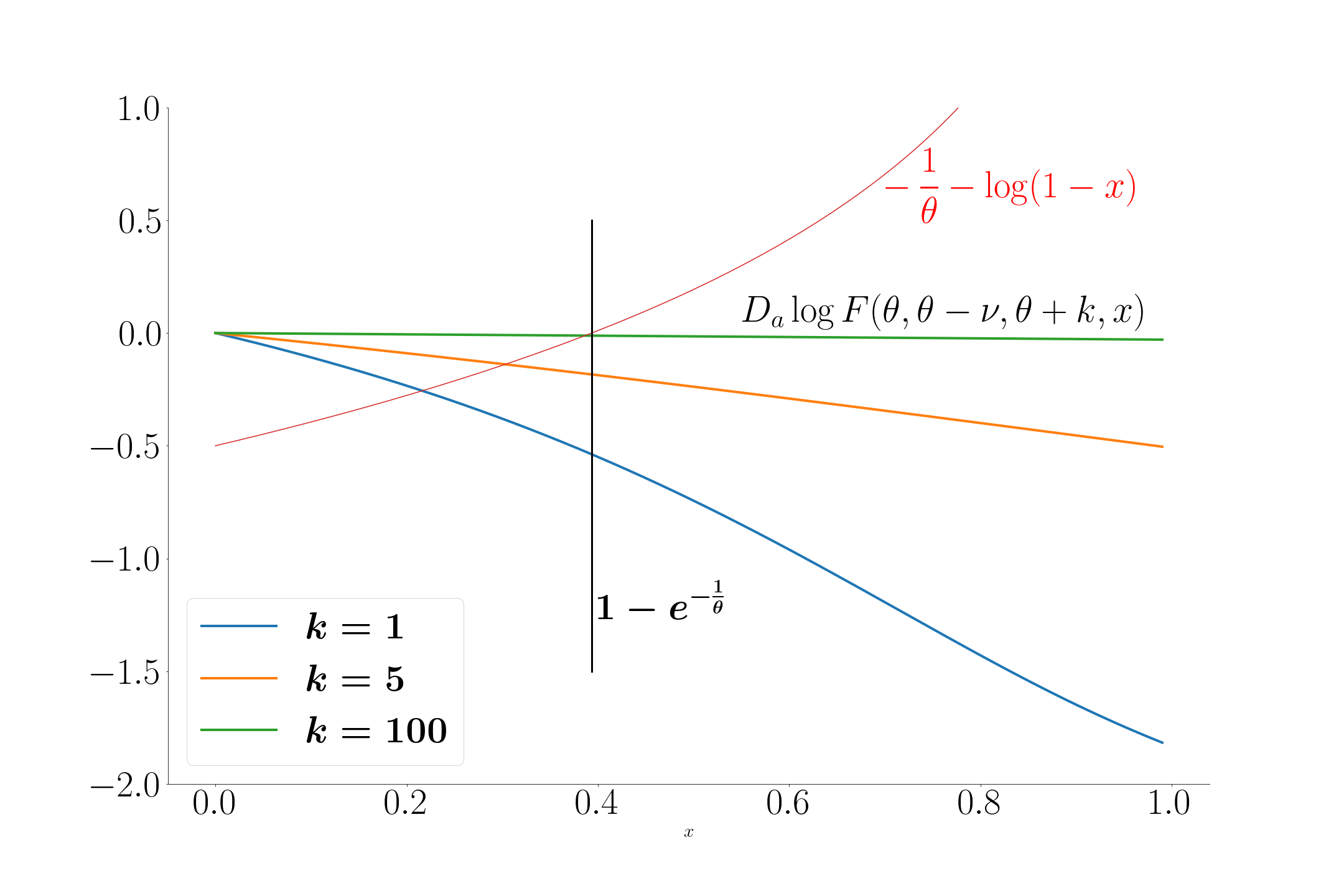}
		\caption{$\nu=5$, $\theta=2$: $\alpha_k \uparrow \alpha^* $ }
		\label{fig:theta 1 nu 3}
	\end{subfigure}
	\hspace*{-1cm}
	\begin{subfigure}[b]{0.5\textwidth}
		\centering
		\includegraphics[width=\textwidth]{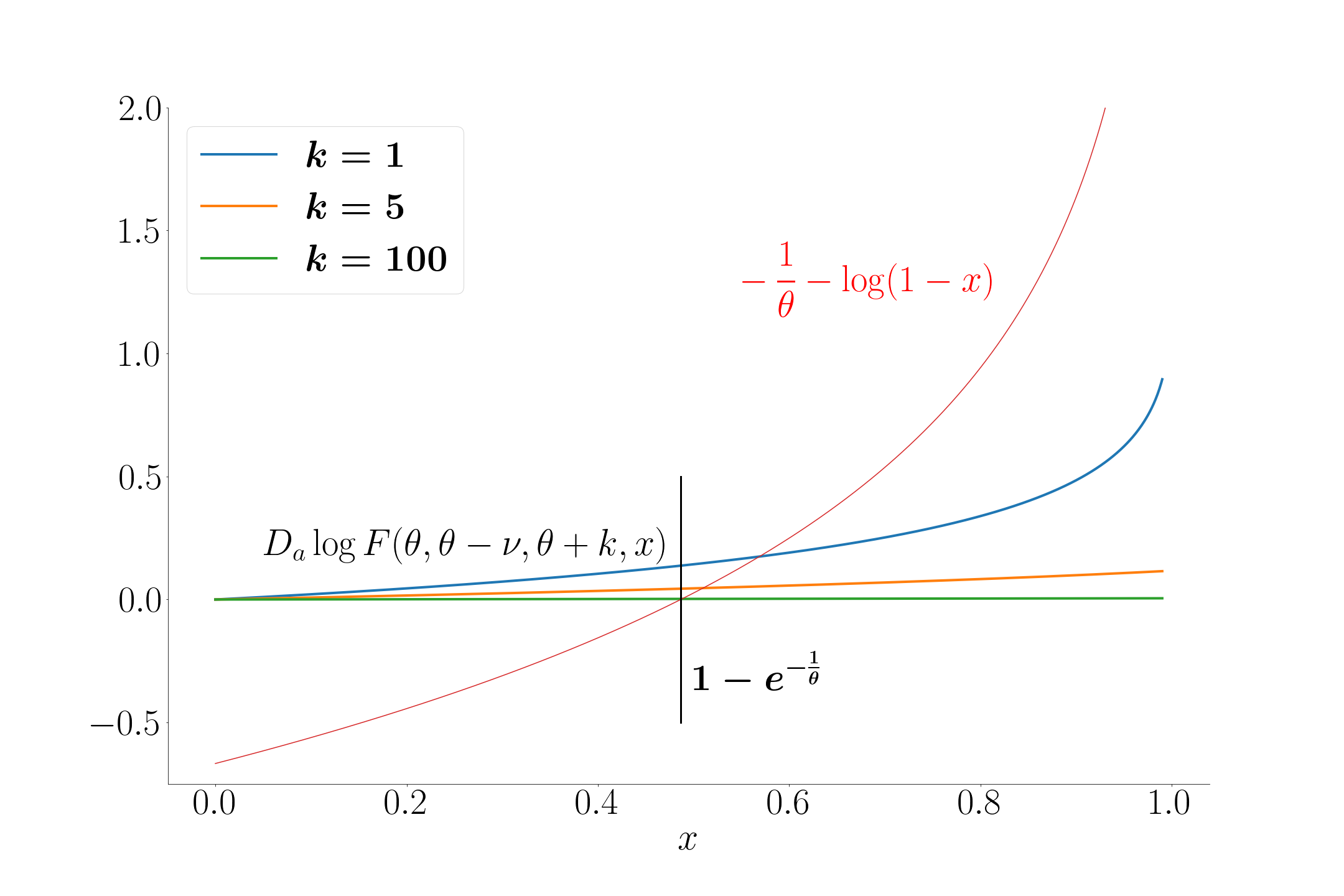}
		\caption{$\nu=1$, $\theta=1.5$: $\alpha_k \downarrow \alpha^*$}
		\label{fig:theta 2 nu 5}
	\end{subfigure}
	
	\caption{Plot of the two sides of equation (\ref{Eq-alpha})}
	\label{fig:Eq-alpha}
\end{figure}

\begin{thm}\label{main}
 The roots satisfy
\begin{eqnarray*}
\alpha_k \uparrow \alpha^* ~~{\rm if~~} \nu&>&\theta,\\
\alpha_k \downarrow \alpha^* ~~{\rm if~~} \nu&<&\theta,\\
\alpha_k = \alpha^* ~~{\rm if~~} \nu&=&\theta.
\end{eqnarray*}
Therefore $\tau^*$ is optimal if $\nu\geq\theta$, and  not optimal if $0<\nu<\theta$.

\end{thm}

\begin{proof}
Straight from the definition of $F$ by the power series (\ref{F}), 
\begin{eqnarray*}
D_a \log (F(\theta,0,\theta+k,x))&=&0,\\ 
D_a \log (F(\theta,\theta-\nu,\theta+k,0))&=&0,\\
\lim_{k\to\infty} D_a \log (F(\theta,\theta-\nu,\theta+k,x))&=&0.
\end{eqnarray*}   
This immediately implies the assertion in the case $\theta=\nu$ and the convergence $\alpha_k\to\alpha^*$ in general.

Suppose $\nu>\theta$. To justify $\alpha_k<\alpha_{k+1}$,  it is sufficient to show that 
$$D_a \log (F(\theta,\theta-\nu,\theta+k,x))<D_a \log (F(\theta,\theta-\nu,\theta+k+1,x)), ~~~x\in(0,1),$$
which by monotonicity of the logarithm is equivalent to 
$$D_a \,\left[\frac{F(\theta,\theta-\nu,\theta+k,x)} {F(\theta,\theta-\nu,\theta+k+1,x)}\right]<0.$$
But this holds  by Lemma \ref{L2-HH}.

The case $\nu>\theta$ is completely analogous.\end{proof}

Since the cutoffs of  $\tau^*$  are related to the roots via (\ref{a-alpha}), we have that $a_k$'s converge to  the value   $a^*:=(1-\alpha^*/q)_+$.
For $\nu\geq \theta$,   the myopic strategy  never stops before time $a^*$.
For $\nu<\theta$,  the {\it optimal} strategy (which is not $\tau^*$) would stop at any success after time $a^*$ regardless of the index.

To contrast this result with another setting, we recall that for the best-choice problem with the Poisson prior $\pi_n=e^{-q} q^k/k!$,
there are only finitely many nonzero cutoffs for each $q>0$, because
the analogues of $\alpha_k$'s  do not accumulate, rather grow about linearly with $k$
(see \cite{BG, BrussCZ, CZ} and Section \ref{UMA}).

\begin{figure}[h]
	\centering
	\begin{subfigure}[b]{0.5\textwidth}
		\centering
		\includegraphics[width=\textwidth]{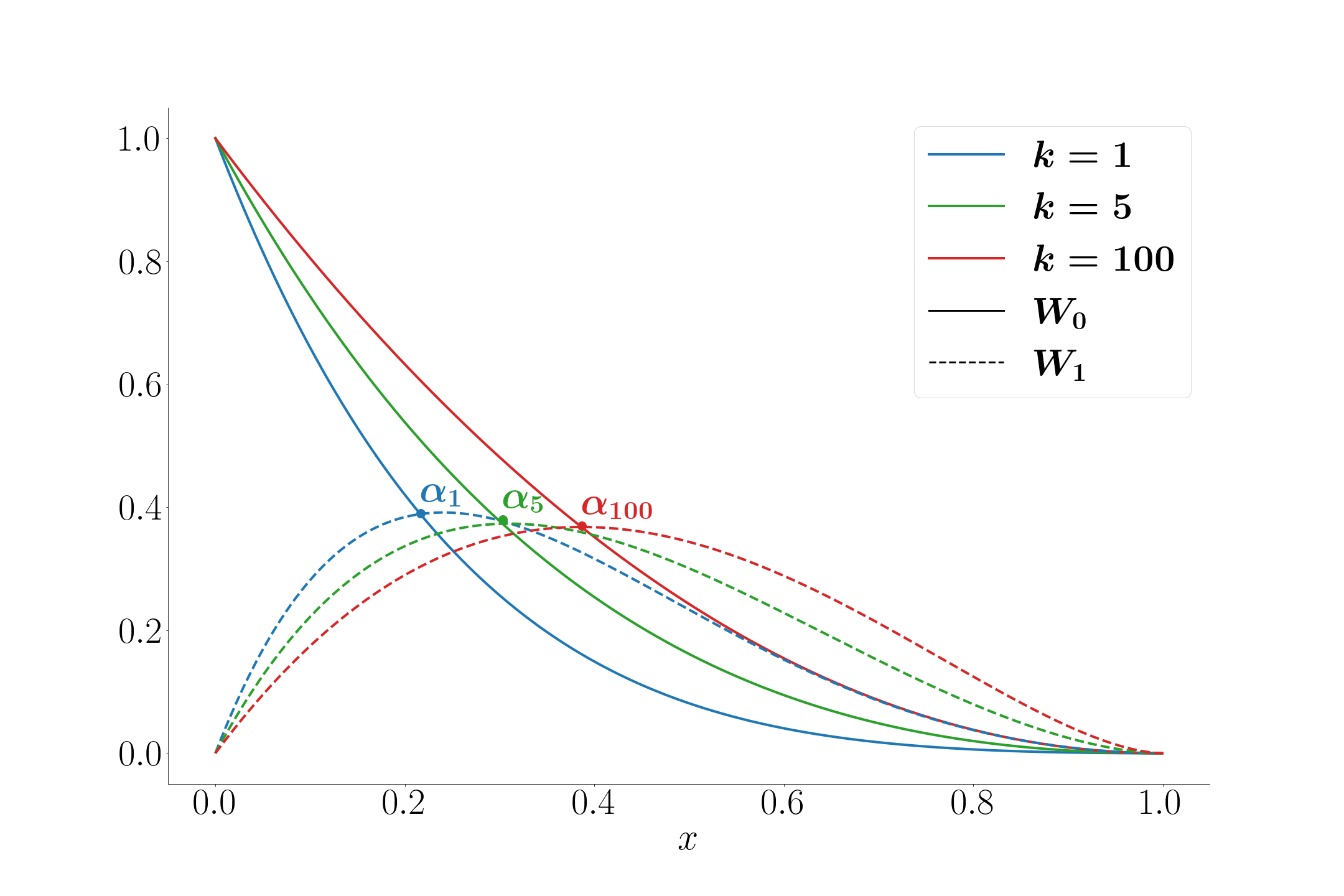}
		\caption{$\nu=5$, $\theta=2$}
		\label{fig:theta 2 nu 5}
	\end{subfigure}
	\hspace*{-1cm}
	\begin{subfigure}[b]{0.5\textwidth}
		\centering
		\includegraphics[width=\textwidth]{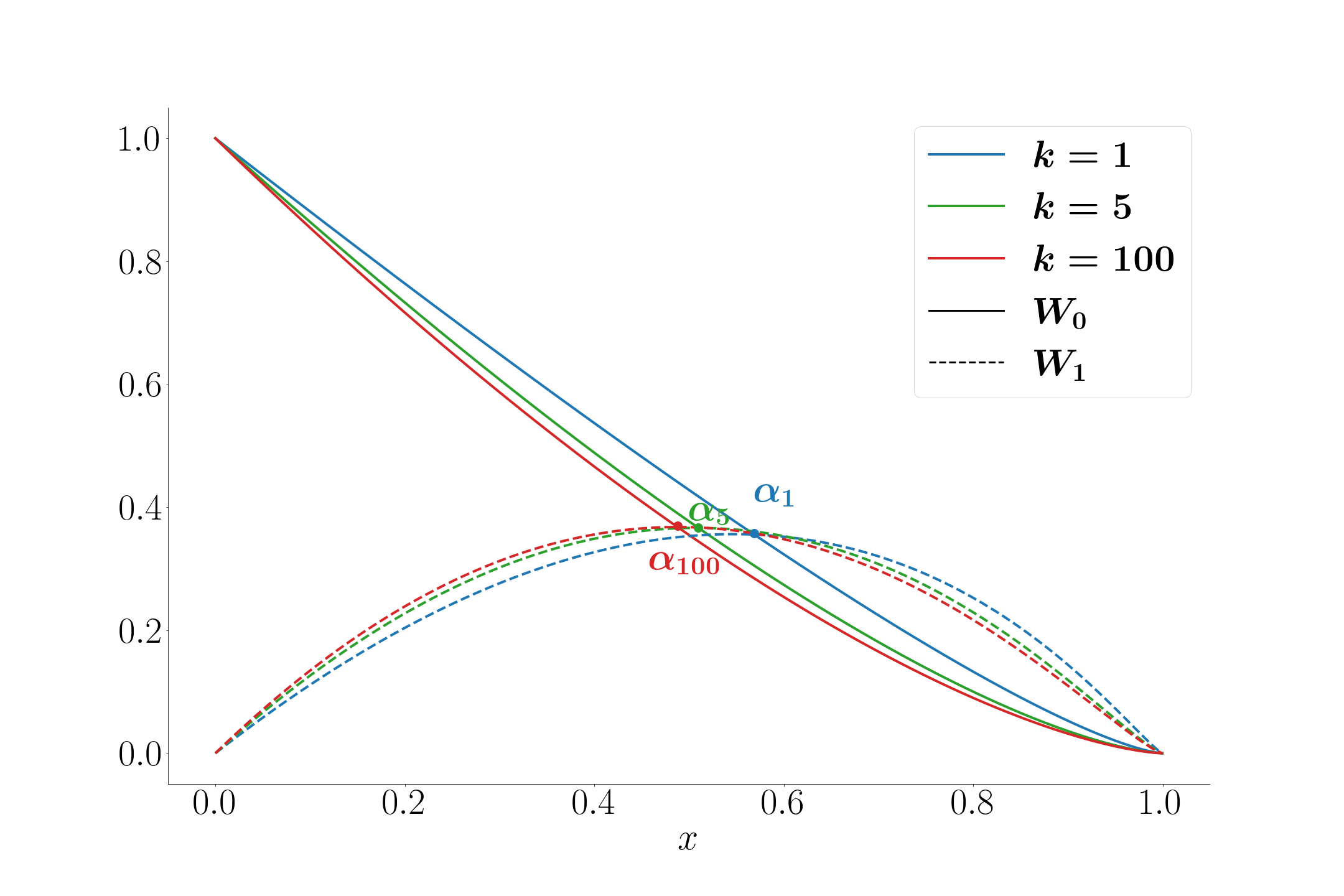}
		\caption{$\nu=1$, $\theta=1.5$}
		\label{fig:theta 1 nu 3}
	\end{subfigure}
	\caption{The winning probabilities $W_0$ and $W_1$}
	\label{fig:S0S1}
\end{figure}

\begin{table}[h]
	
	\setlength{\tabcolsep}{5.8mm}

	\centering
	\begin{tabular}{lcc}\toprule
		
		\multirow{2}{*}[-3pt]{ \boldmath{$k$}}& \boldmath{$\nu=5, \theta=2$}& \boldmath{$\nu=1, \theta=1.5$ }\\

		\cmidrule(lr){2-2}\cmidrule(lr){3-3}
		& {$\alpha_k \uparrow \alpha^*=0.393469$ }  & {$\alpha_k \downarrow \alpha^*=0.486583$}  \\
		\midrule
		
		1 &  0.216390 & 0.568837 \\
		2 & 0.249979 & 0.537613 \\
		3 & 0.273297 & 0.523253 \\
		4 & 0.290259 & 0.515103 \\
		5 & 0.303095 &  0.50988 \\
		6 & 0.313126 & 0.506259 \\
		7 & 0.321171 & 0.503604 \\
		8 & 0.327761 & 0.501577 \\
		9 & 0.333257 & 0.499979 \\
		10 & 0.337909 & 0.498687 \\
		20 & 0.362155 & 0.492737 \\
		50 & 0.379922 & 0.489067 \\
		100 & 0.386508 & 0.487829 \\
		1000 & 0.392755 & 0.486708 \\
		10000 & 0.393398 & 0.486595 \\
		100000 & 0.393462 & 0.486584 \\
		$\infty$ &  0.393469&  0.486583\\
		\bottomrule
	\end{tabular}
	\caption{Critical points {\rm $\alpha_k$:} Solution to $W_1(x, k)=W_0(x, k)$}
	%
			%
			%
			%
			%
	
	\label{tab:roots}
\end{table}

\paragraph{The case $\nu=\theta$: a  Poisson process of success epochs.}
In the case $\nu=\theta$ there is much simplification, as (\ref{PGF1}) for every $k\geq0$ becomes
$z\mapsto (1-x)^{\theta(1-z)},$
which is the  probability generating function of the Poisson distribution with mean $-\theta\log(1-x)$.
All roots  $\alpha_k$ in this case coincide with $\alpha^*$, hence   $\tau^*$ is a single cutoff strategy which skips
all trials before
$a^*$
then stops at the first success. 
The optimal winning probability is 
$${\cal S}_1(a^*,k)=-\theta (1-x)^\theta\log(1-x), ~~~{\rm for~~}x=q-(q-\alpha^*)_+,$$
hence equal to $1/e$ for $q$ sufficiently large (more precisely, in the range $1-e^{-1/\theta} \leq q<1$). 
This generalises the $\theta=1$ result from \cite{BrussSam1} for the  geometric prior.

A probabilistic reason for such simple solution is the following.
From 
Proposition \ref{Structure} (ii) and (\ref{profile}),  the rate function for the point process of success epochs does not depend on index $k$.
$$   \frac{k+\theta}{t+q^{-1}-1} p_{k+1}=   \frac{\theta}{t+q^{-1}-1},$$

Thus,  the counting process of success epochs has deterministic infinitesimal compensator, therefore the process is Poisson 
by  Watanabe's martingale characterisation. 
So  the last-success problem in this case is equivalent to the elementary problem of maximising the probability of stopping at the last epoch of a Poisson process.


\paragraph{The  polynomial case: $\nu-\theta\in{\mathbb N}$.}
If $\nu-\theta=m$ for some positive integer $m$, the hypergeometric series (\ref{PGF1}) is a  polynomial in $x$ 
of degree $m$. Similarly to Jacobi polynomials ($P^{\alpha,\beta}$ with positive parameters) the function is
representable via Rodrigues' formula as
$$
F(\theta-\theta z,-m,\theta+k,x)= \frac{x^{1-\theta-k}(1-x)^{k+\theta z+m}}{(\theta+k)_m} D_x^m \left[\frac{x^{\theta+k+m-1}}{(1-x)^{k+\theta z}} \right],
$$
where $D^m_x$ denotes the $m$th  derivative in $x$.

Showing that  the roots increase in this polynomial case can be concluded by a shorter argument based on 
 the parameter derivative formula due to Fr{\"o}hlich
(see Lemma \ref{L3-HH} below), 
$$D_a\log F(\theta,-m,\theta+k,x)= \sum_{j=0}^{m-1}\frac{1}{\theta+j} -\sum_{j=0}^{m-1} \frac{m!}{j!(m-j)} \frac{F(\theta,-j,\theta+k,x)}{F(\theta,-m,\theta+k,x)},$$
and applying the  third part of Lemma \ref{L1-HH} to check 
 that the right-hand side is increasing in $k$.

The left-hand side of  equation  (\ref{Eq-alpha}) for $\alpha_k$ becomes increasingly more complicated with growing degree $m=\nu-\theta$. To illustrate:

$$D_a\log F(\theta,-m,\theta+k,x)= \begin{cases}
    \dfrac{x}{    \theta x- \theta-k }\,,{\rm ~~for~}m=2,\\
 \dfrac{(2\theta+1)x^2-2(\theta+k+1)x}{\theta(\theta+1)x^2- 2\theta (\theta+k+1) x+(\theta+k)(\theta+k+1)}\,,{\rm ~ ~for~}m=3.
\end{cases}
$$


\paragraph{The case $\nu=0$.} The edge case $\nu=0$ corresponds to the logarithmic series prior
\begin{equation}\label{LSP}
\pi_n=   \frac{q^n}{ |\log(1-q)| n}\,, ~~~n\geq1,
\end{equation}
with $0<q<1$. Formula (\ref{PGF1}) is still valid, along with the argument for $\alpha_k\downarrow\alpha^*$, hence $\tau^*$ is not optimal.
The instance $\theta=1$ was studied in much detail in \cite{GD}.

\section{The winning probability}

We first assume arbitrary success profile  $\profile$ and prior $\prior$.
Let $\tau$ be a strategy with  a 
stopping set
$B=\bigcup\limits_{k=1}^\infty ([b_k,1)\times \{k\}),$
where the cutoffs satisfy $1> b_1\geq b_2\cdots\geq0$.
By the monotonicity, once the process $(t, N_t)$ enters $B$ it stays there all the remaining time.
We define the {\it precursor} associated with $\tau$  to be the random time
\begin{equation}\label{sigma1}
\sigma:= 
\min\{t: (t, N_t+1)\in B\}.
\end{equation}
Clearly, $\sigma<\tau$ and  $\sigma\leq b_1<1$. 
The rationale behind this definition is that starting from time $\sigma$, the strategy $\tau$ acts  
greedily, by stopping 
 at the first available  success.  
Hence the winning probability with $\tau$ can be expressed as ${\mathbb E} {\cal S}_1(\sigma, N_\sigma)$.

Both $\tau$ and
  $\sigma$  are stopping times adapted to the filtration $({\cal F}_t, ~t\in[0,1])$ but, as the next lemma shows, 
$\sigma$ takes values in a larger set $\{T_1,\cdots,T_N, b_1, b_2,\cdots\}$.

\begin{lemma}\label{L-sigma}     $N_{\sigma}=k$ for  $k\geq1$ holds iff either
\begin{itemize}
\item[\rm(i)] $N_{b_{k+1}}=k$, in which case $\sigma=b_{k+1}$,
\end{itemize}
or
\begin{itemize}
\item[\rm(ii)] $N_{b_{k+1}}<k\leq N_{b_k}$, in which case $b_{k+1}<\sigma=T_k\leq b_k$.
\end{itemize}
Also, $N_\sigma=0$ is equivalent to $N_{b_1}=0$.  
\end{lemma}
\begin{proof}

In the case (i) we have $(b_{k+1}, k+1)\in B$, and for $t<b_{k+1}$ it holds that  $N_t\leq k$ hence $(t, N_t)\not\in B$.
Thus $(\sigma, N_\sigma)=(b_{k+1},k)$.

In the case (ii) we have $b_{k+1}<T_k\leq b_k$, so $(T_k, k+1)\in B$. On the other hand, for $t<T_k$  we have  $N_t\leq k-1, t<b_k$ hence  $(t,N_t)\not\in B$.
Thus $(\sigma, N_\sigma)=(T_k,k)$.

Note that  the condition `(i) or (ii)'  means that $N_t=k$ for some $t\in [b_{k+1}, b_k]$. This condition is also necessary for $N_\sigma=k$. 
Indeed, if $N_{b_{k+1}}>k$ then $T_{k+1}\leq b_{k+1}$ together with $N_\sigma=k$ would imply $\sigma<T_{k+1}\leq b_{k+1}$ and $(\sigma, N_\sigma+1)\notin B$;
a contradiction with (\ref{sigma1}). 
Likewise, if $N_{b_k}<k$ then $N_\sigma=k$ would imply $b_k<T_k\leq \sigma$, hence choosing $t$ in the range $\max(b_k, T_{k-1})<   t<T_k$ we have 
$(t, N_t+1)=(t, k)\in B$, which in the view of  $t<\sigma$ contradicts the minimal property in  (\ref{sigma1}).

In the case $k=0$, only option (i) can occur.
\end{proof}

\paragraph{Remark} 
On the event $N=n$ the precursor is
\begin{equation*}\label{sigma2}
\sigma= \min\{t:~   N_t+ \#\{j\geq 1:~ b_{n-j+1}\leq t\}=n\}.
\end{equation*}
This has a nice combinatorial interpretation.
 The cutoffs $ b_n\leq b_{n-1}\leq\cdots\leq  b_1$ split $[0,1]$ into $n+1$ compartments, where the 
points $T_1,\cdots, T_n$ are allocated. 
The sequence of $n+1$ cluster sizes can be written in stars and bars notation, like for instance
$ *|~|**{\bf\color{red}|}~|*|**$
 for $n=6$. 
The precursor corresponds to  the  $n^\text{th}$ symbol (star or bar) in this sequence.

The winning probability can be computed by conditioning on the number of trials as
\begin{eqnarray}\label{WinP}
{\mathbb E} {\cal S}_1(\sigma, N_\sigma)= \sum_{n=1}^{\infty} \pi_n  \sum_{k=0}^{n-1} {\mathbb P}(N_\sigma=k|N=n) s_1(k+1,n).
\end{eqnarray}
The formula is valid for arbitrary $\profile$ and $\prior$.
To apply (\ref{WinP})  one needs to find the distribution of $N_\sigma$ given $N=n$ trials.

\begin{lemma} For $0\leq k\leq n$, we have 
\begin{eqnarray*}
{\mathbb P}(N_\sigma=k|N=n)&=&\\
{n\choose k} b_{k}^k (1-b_{k})^{n-k} &+& \sum_{i=k+1}^n {n\choose i} \left[ b_k^i(1-b_k)^{n-i}-b_{k+1}^i (1-b_{k+1})^{n-i}\right], 
\end{eqnarray*}
where $b_0:=1$.
\end{lemma}

\begin{proof} The case $k=0$ is easy:
${\mathbb P}(N_\sigma=0\,|\,N=n)=(1-b_1)^n$. For $k\geq 1$, 
by Lemma \ref{L-sigma}, the event occurs if  
either $N_{b_{k+1}}=k$ or   $T_k$ falls in $[b_{k+1}, b_k]$, hence
\begin{eqnarray*}
{\mathbb P}(N_\sigma=k|N=n)&=&
{\mathbb P}(N_{b_{k+1}}=k|N=n)+ {\mathbb P} (N_{b_{k+1}}<k \leq N_{b_{k}}|N=n)=\\
{\mathbb P}(N_{b_{k}}=k|N=n)&+& {\mathbb P} (N_{b_{k+1}}<k+1 \leq N_{b_{k}}|N=n)=\\
{\mathbb P}(N_{b_{k}}=k|N=n)&+&{\mathbb P} (N_{b_{k}}\geq k+1|N=n)- {\mathbb P} (N_{b_{k+1}}\geq k+1|N=n),
\end{eqnarray*}
and the asserted formula  follows.
\end{proof}

Formula (\ref{WinP}) resembles known formulas for the winning probability in the full-information best-choice problem \cite{GM, Poros1}.
Alternatively, a direct application of Lemma \ref{L-sigma} yields

\begin{equation}\label{WinP2}
{\mathbb E} {\cal S}_1(\sigma, N_\sigma)=\sum_{k=0}^\infty {\mathbb P}(N_{b_{k+1}}=k){\cal S}_1(b_{k+1},k)+
\sum_{k=1}^\infty \int_{b_{k+1}}^{b_k} {\cal S}_1(t,k){\mathbb P}(T_k\in {\rm d}t).
\end{equation}

\begin{figure}[h]
	\centering
	\includegraphics[width=10 cm]{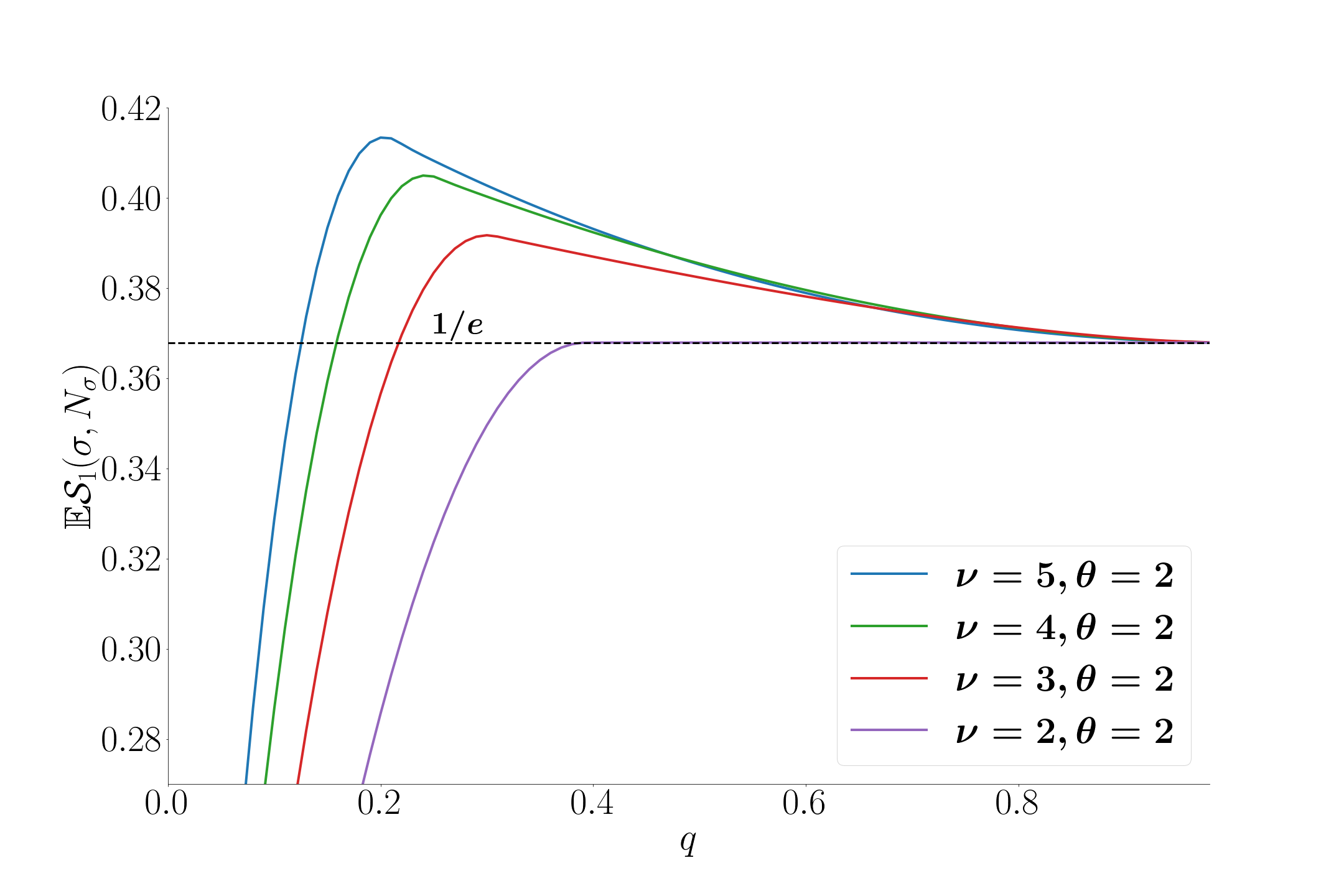}
	\caption{Winning probability as a function of $q$ for different configurations of  $\nu$ and $\theta$}%
	\label{F1}
\end{figure}

In  the case of ${\rm NB}(\nu,q)$ prior, the distribution of $N_t$ is given by
Proposition \ref{Structure},  which also yields the density of $T_k$
$$
\prob(T_k\in {\rm d}t)=\prob(N_t=k-1)\frac{k+\nu-1}{t+q^{-1}-1}{\rm d}t
$$
needed to apply (\ref{WinP2}).

\section{Asymptotics and limit forms of the problem}

\subsection{Asymptotic optimality} 
For profile (\ref{profile}), there is a universal   asymptotics valid under the sole assumption that the number of trials becomes large.


\begin{assertion} As $N\to \infty$ (in probability), the strategy with single cutoff $e^{-1/\theta}$ stops at the last success with 
probability approaching $e^{-1}$, which is the  asymptotically optimal winning probability.
\end{assertion}
\begin{proof} Condition on $N=n$, with $n$ large. By  the second limit relation in (\ref{fixed-n}) the  winning probability cannot exceed  $e^{-1}+\varepsilon$,
even if  the observer knows  $n$  apriori. Hence $e^{-1}$ is an asymptotic upper bound.

On the other hand, it is still true that $s_1(K+1,N)\to e^{-1}$ for $K/N$ converging in probability to $e^{-1/\theta}$, as $N\to\infty$.
By the law of large numbers, 
this holds for
$K=N_{e^{-1/\theta}}$,   the count of trials before time $e^{-1/\theta}$.
 Therefore  the strategy with single cutoff ${e^{-1/\theta}}$ has the winning probability converging to 
$e^{-1}$.
\end{proof}

In fact, there are many Markovian asymptotically optimal strategies.

\begin{assertion}\label{many} Let $\tau$ be a strategy with cutoffs $b_1, b_2,\cdots$. If $~b_k>0$ and $\lim_{k\to\infty}b_k= e^{-1/\theta}$ then $\tau$ is asymptotically optimal as $N\to\infty$ (in probability).
\end{assertion}

\begin{proof}  Let $\tau'$ be the strategy with cutoff $e^{-1/\theta}$.
We argue that the equality $\tau=\tau'$  holds with probability approaching one.
Indeed, the least cutoff  $b:=\min b_i$  is strictly positive.
Fix  small $\varepsilon$ and choose  $k$ such that $|b_j-e^{-1/\theta}|<\epsilon$ for $j>k$.   If $N_b>k$ 
 then $\tau=\tau'$ unless there is a success epoch $(T_i, i)^\circ$ in the $\varepsilon$-vicinity of $e^{-1/\theta}$.   
As $N\to\infty$ also $N_b\to\infty$ and the probability that there exists such $T_j$ is of the order $O(\varepsilon)$ (see the next proposition).
\end{proof}

\subsection{A Poisson limit problem}\label{PL}
More insight in the asymptotics is gained by looking at the temporal pattern of success trials.
Let  $\Pi$ be a nonhomogeneous Poisson point process with  rate function $\theta/t$, $t\in (0,1]$.
This is the familiar selfsimilar Poisson process \cite{PitmanYou} restricted to the unit interval.
Note that $\Pi$ has infinitely many points and these  accumulate near $0$.

\begin{assertion}\label{assert1} As $N\to \infty$ (in probability), the
point process of success epochs converges weakly to $\Pi$.
\end{assertion}
\begin{proof} 
Condition on $N=n$,  for $n$ large. The approximation of the point process of success epochs  by $\Pi$ works  for trials occurring at fixed times $k/n, 1\leq k\leq n,$ 
 since $\sum_{k=1}^\infty p_k=\infty$ and $\sum_{k=1}^\infty p_k^2<\infty$, see \cite{Pfeifer}.
This is readily extended to the mixed binomial process  by noting that  by the law of large numbers  the process $(N_t/N, \, t\in[0,1])$ converges weakly to the identity function.
\end{proof}

 The following last-arrival problem can be seen as a limit form of the last-sucess problem.
Suppose the points of $\Pi$ are revealed in their natural order with the objective to stop at the 
last point before time $1$.  The optimal strategy is then to stop at the first point after time $e^{-1/\theta}$.
The optimality is justified by a standard monotone case argument, with the cutoff found by equating the probability $t^\theta$ of zero Poisson points on $[t,1]$
and the probability $-(\theta \log t  )     t^\theta$
of exactly one such point.
For $\theta=1$, this limit form is a special case of the  
Gianini and Samuels  \cite{GianiniSam} infinite secretary problem.

\subsection{Infinite priors} 
By the virtue  of  Proposition \ref{Structure},
the last-success problem with negative binomial prior can be formulated straight in terms of the P{\'o}lya-Lundberg process 
starting from  $N_0=0$, that is with  the initial state $(0,0)$. 
The optimal strategy is characterised by some stopping set of states
$B$. If we consider the pacing process 
starting from any other initial state $(t_0,k_0)\in(0,1)\times\{0,1,\cdots\}$, the optimal strategy will be
  to stop at the first subsequent success $(t,k)^\circ$ with $(t,k)\in B$.
In this sense the optimal strategy does not depend on the initial state.
This follows from the optimality principle and the fact that every state is reachable from $(0,0)$.
Likewise, the myopic strategy is determined by the stopping set 
 $A$ as in (\ref{A}) regardless of the initial state of the pacing process.

Let  $\tau^\dagger$ be the strategy with  cutoffs $1-\alpha_k, k\geq1$. 
By Proposition \ref{many}   this strategy    is asymptotically optimal. In particular,  $\tau^\dagger$ is asymptotically optimal under the negative binomial prior 
${\rm NB}(\nu,q)$ as $q\to 1$. We aim to confer to $\tau^\dagger$ yet another sense of optimality in the case $\nu\geq\theta$.

To that end, note that for  initial state $(t_0, k_0)$ with $t_0\neq0$ the P{\'o}lya-Lundberg process with 
rates formula (\ref{JR}) is well defined also for $q=1$. The case $t_0=0$ must be excluded  because the rate function has a pole at $0$.
In this setting, $\tau^\dagger$ is the myopic strategy whose winning probability depends on $(t_0,k_0)$ and converges to $e^{-1}$ as $t_0\to0$.
If $\nu\geq \theta$, strategy $\tau^\dagger$  is optimal in the last-success problem with any such initial state.

The $q=1$ setting can be interpreted as  the last-success problem with infinite prior 
\begin{equation}\label{InP}
\pi_n=\frac{(\nu)_n}{n!}, ~n\geq0.
\end{equation}
Although this is not a probability distribution,  a formal Bayes calculation of the posterior  for $N-N_{t_0}$  given $N_{t_0}=k_0$ (with $t_0\neq 0$) yields the proper negative binomial 
distribution ${\rm NB}(\nu+k_0, 1-t_0)$, in accord with  part (ii) of 
Proposition \ref{Structure}. Exactly this  distribution of $N-N_{t_0}$   emerges  under the P{\'o}lya-Lundberg process with parameter $q=1$ and initial state $(t_0,k_0)$.

For the best-choice problem ($\theta=1$)  with $\nu=1$
the idea of formal Bayes optimality  was introduced by Stewart \cite{Stewart}, who proved the optimality of    the $1/e$-strategy (which coincides with $\tau^\dagger$ in this case).
For $\nu=1$,  (\ref{InP}) is the uniform distribution on nonnegative integers, regarded  in the statistical literature as `noninformative prior'
to address the situation of complete ignorance about the value of $N$.
Bruss \cite{CZ} considered the best-choice problem with mixed Poisson pacing process whose  rate was assumed  to follow the infinite uniform distribution on positive reals, however 
it is not hard to see that
this setting is equivalent to that of \cite{Stewart}.
Thus the $1/e$-strategy is both Bayes optimal under the infinite uniform prior and optimal in the limit model of Section \ref{PL}
(see discussion in 
  \cite{BrussSam1}, p. 829).
But note that in our  two parameter  framework the latter applies  to $\tau^\dagger$ only if $\theta=\nu$.

The edge case $\nu\to0$ of (\ref{InP}) corresponds to the infinite logarithmic series prior $\pi_n=n^{-1},~n\geq0$, which is  the $q=1$ analogue of (\ref{LSP}). 
In different terms, the model 
was proposed in \cite{BrussYor} as yet another way to model the complete ignorance about $N$. 
Bruss and Rogers \cite{BrussRogers21} observed that the pacing process is a birth process and showed that 
the $1/e$-strategy is not  optimal in the best-choice problem with initial state $(t_0,1)$, $t_0<1/e$.
To this we can add 
 that the monotone case of optimal stopping fails for every $\theta>0$, thus neither myopic nor any single cutoff strategy is optimal
for initial states $(t_0,k), t_0<e^{-1/\theta}$.

\section{A unified Markovian approach}\label{UMA}

A simple time change allows us  to embed the last-success problems with different $q$  in a unified framework of a single Markov process.
Parameter $q$ will assume the role of the initial state.
Consider the general  power series prior (\ref{p-series}), and let $I\subset[0,\infty)$ be the maximal convergence interval 
of  $\sum_n w_n q^n$.
If  $I=[0,\rho)$  for some $\rho<\infty$ then there exists an  improper prior $\pi_n=w_n\rho^n,~n\geq 0$. 
The latter does not apply, for instance, to the shapes like  $w_n=1/n!$ or $w_n=1/n^2$.

The associated
 mixed binomial process $(N_t, ~t\in[0,1])$ has the property that
the conditional distribution of $N$ given $N_t=k$ depends on $q$ and $t$ through 
$x=q(1-t)$. Thus the process
$(N_{1-x/q}, ~q\geq x\geq0)$, with the temporal variable $x$ decreasing from $q$ to $0$,  has 
 a jump rate $r(x,k)$ not depending on $q$.
Let  $M$  be  a Markov process with the state space $S=I\times\{0,1,\cdots\}$ and the transition kernel defined by 
$${\mathbb P}(M_{x-{\rm d}x}=(x-{\rm d}x,k+1)\,| M_x=(x,k))=r(x,k){\rm d}x.$$
The process $M$ starting from $(q,0)$ has the same distribution as $(N_{1-x/q}, ~q\geq x\geq 0)$.
Suppose state $(x,k)\in S$ of  $M$ is marked  as success with probability $p_k$, independently of anything else. 
The last-sucess problem with parameter $q$ translates as the optimal stopping of $M$ starting from state $(q,0)$.

A Markovian stopping strategy in this setting is determined by  a subset of  $S$.
Denote $V(x,k)$ the maximum winning probability achievable by proceeding  from state $(x,k)$, $W_0(x,k)$ the probability of no successes following the state  $(x,k)$ and
let
$C:=\{(x,k)\in S: W_0(x,k)\geq V(x,k)\}$. The optimal strategy by starting from any given $(x,k)\in S$ is to stop at the first state in $C$ which is marked as success. 
This very set $C$ corresponds to the Bayes optimal strategy whenever the infinite prior exists.

Specialising to  the negative binomial prior (\ref{NB}) (with $\nu$ fixed) we have $I=[0,1)$,
 $r(x,k)=(k+\nu)/(1-x)$,
\begin{equation}\label{C}
C=\bigcup_{k=1}^\infty ([0,\alpha_k]\times\{k\})
\end{equation}
Similarly to \cite{GD}, the value function satisfies the optimality equation
$$D_x V(x,k)=\frac{k+\nu}{1-x} \{V(x,k+1)-V(x,k)\}+\frac{(k+\nu)\theta}{(\theta+k)(1-x)}\{W_0(x,k+1)-V(x,k+1)\}_+$$
with the initial condition $V(0,k)=0$. Conditioning on the time of the first subsequent trial yields an equivalent integro-differential equation

\begin{eqnarray*}
V(x,k)=~~~~~~~~~~~~~~~~~~~~~~~~~~~~~~~~~~~~~~~~~~~~~~~~~~~~~~~~~~~~~~~~~~~~~~~~~~~~\\
 \int\limits_0^1 \left[ \frac{\theta}{\theta+k} \max(W_0(y,k+1),V(y,k+1))+\frac{k}{\theta+k} V(y,k+1)     \right]   
 \frac{(\nu+k)x(1-x)^{\nu+k}}{(1-xy)^{\nu+k+1}} {\rm d}y.
\end{eqnarray*}
The role of the boundary condition at `$k=\infty$' is played by the asymptotics
$$\lim_{k\to \infty} V(k,x)=\begin{cases} ~~~ e^{-1},~~~{\rm for}~~~~1-e^{-1/\theta}\leq x\leq1,\\
-\theta(1-x)^\theta\log(1-x),~{\rm for}~0\leq x\leq 1-e^{-1/\theta},  \end{cases}
$$
which follows from Proposition \ref{assert1}.

\begin{figure}[ht]
	\centering
	\begin{subfigure}[b]{0.5\textwidth}
		\centering
		\includegraphics[width=\textwidth]{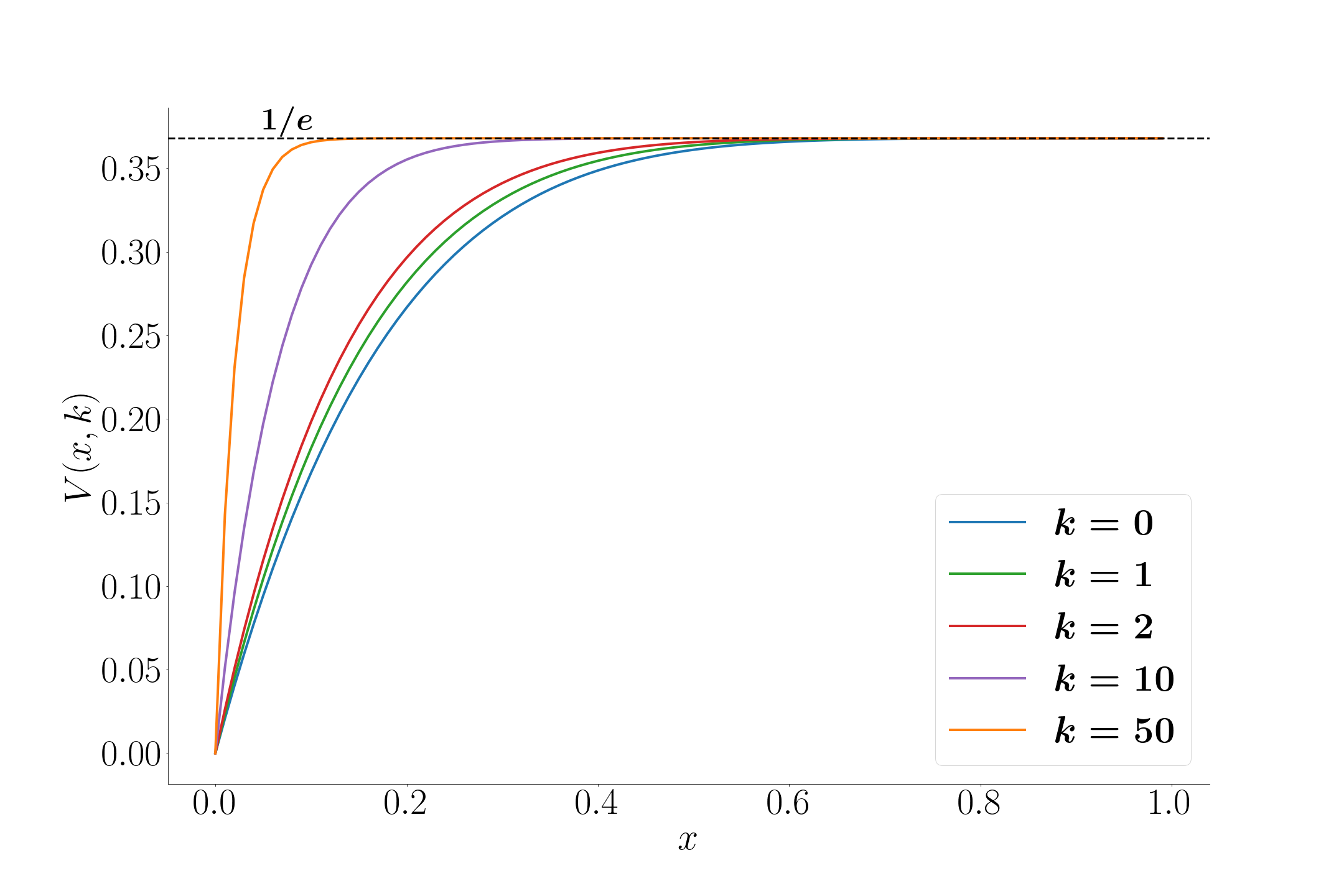}
		\caption{$\nu=5$ , $\theta=2$}
		\label{fig:theta 1 nu 3}
	\end{subfigure}
	\hspace*{-1cm}
	\begin{subfigure}[b]{0.5\textwidth}
		\centering
		\includegraphics[width=\textwidth]{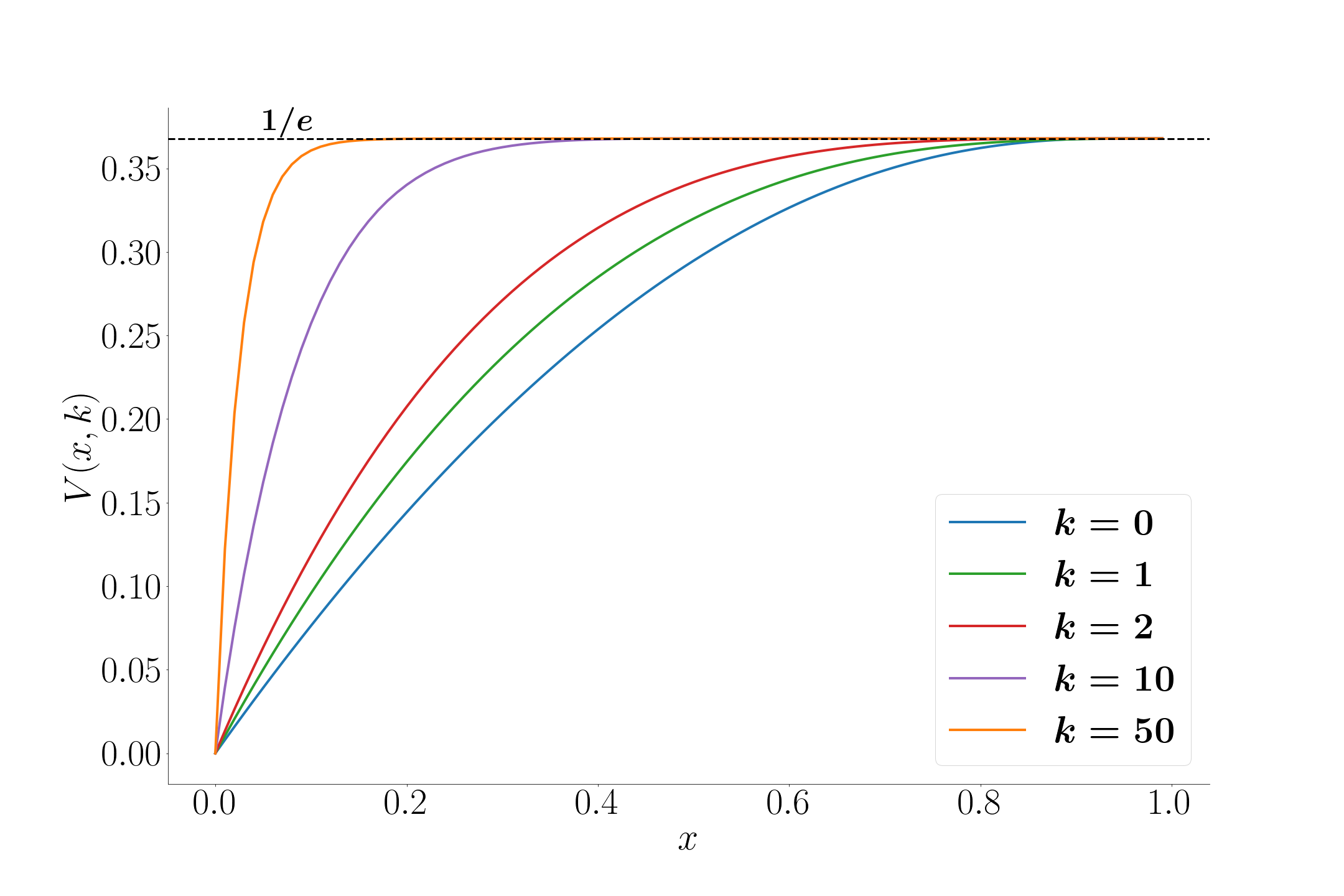}
		\caption{$\nu=1.5$, $\theta=1$}
		\label{fig:theta 2 nu 5}
	\end{subfigure}
	
	\caption{Value Function $V(x, k)$ }
	\label{fig:Eq-alpha}
\end{figure}

Some conclusions about the value function can be drawn 
from what we know about the optimal strategy. 
For $\nu\leq \theta$, we have $\alpha_k\geq \alpha^*=1-e^{-1/\theta}$ hence the optimal strategy stops greedily for $0\leq x\leq \alpha^*$ and in this range we have
$V(x,k)=W_1(x,k), ~~k\geq0$. 
 For $\nu\geq\theta$, 
from  Theorem \ref{main} we have
$V(x,k)= W_1(x,k)$ if $x\leq \alpha_{k+1}, ~k\geq 0$; and since it is never optimal to stop for $x>\alpha^*$ we have a relation
$$V(x,k)=\sum_{j=0}^\infty    \lambda_jV(\alpha^*,k+j), ~~~\alpha^*<x\leq 1, ~k\geq 0,$$
where $\lambda_j$'s are the masses of ${\rm NB}(\nu+k, xe^{-1/\theta}/(1-x+xe^{-1/\theta}))$.

\paragraph{Remark}
To revisit another known example, suppose $\theta=1$ and
$\pi_n= e^{-q}{q^n}/{n!}$, so $(N_t, t\in [0,1])$ 
is the Poisson process with rate $q$. In this case the rate function for $M$ is constant $r(x,k)=1$.
From \cite{CZ} we know that
 the optimal stopping set is of the form (\ref{C}), where $\alpha_k$ is the positive root of 
the equation
$$\frac{1}{k}+\sum_{j=1}^\infty \frac{x^j}{j!(k+j)}\left(1-\sum_{i=1}^j \frac{1}{i+k-1}  \right)=0,$$
 and that
the sequence of roots is strictly increasing. Therefore the optimal strategy coincides with the myopic strategy.
Unlike the negative binomial case, the roots do not converge, rather have 
asymptotics $\alpha_k\sim k(e-1)$ which was observed in \cite{Ciesielski}.
The convergence interval of $\sum_n q^n/n!$ is $I=[0,\infty)$ and the
 `infinite Poisson prior' does not exist.

\section{Appendix}

\subsection{The hypergeometric function}\label{HGF}

This section collects    monotonicity  properties of the Gaussian hypergeometric series
\begin{equation}\label{F}
F(a,b,c,x):= \sum_{j=0}^\infty \frac{(a)_j (b)_j}{(c)_j}  \frac{x^j}{j!}\,,
\end{equation}
viewed as a function of  real parameters $a, b, c$  and $x\in[0,1]$.
We shall use standard formulas
\begin{eqnarray}
F(a,b, c, x)&=&F(b ,a , c, x),\label{symmetry} \\
F(a,b, c, x)&=&(1-x)^{c-a-b} F(c-a, c-b, c, x),\label{change-par}\\
D_x F(a,b,c,x)&=&\frac{ab}{c}\,F(a+1,b+1,c+1,x)\label{der}
\end{eqnarray}
and  Euler's integral representation 
\begin{equation}\label{Euler}
F(a,b, c, x) =\frac{\Gamma(c)}{\Gamma(b) \Gamma(c-b)} \int_0^1 \frac{z^{b-1} (1-z)^{c-b-1}}{(1-xz)^{a}}\,{\rm d}z ~~~{\rm for}~~c>b>0.
\end{equation}

\begin{lemma}\label{L0-HH} For positive $c>a>0$ and $0<x<1$ 
$$
\sgn D_c F(a,b,c,x) = - \sgn b.
$$
\end{lemma}

\begin{proof}
For $b>0$ the assertion is obvious from the monotonicity of $(c)_j$ and the series formula of the hypergeometric function.

For negative $b$ we argue by induction.
First suppose $-1\leq b<0$. Since $b+1>0$,   by the above $F(a+1,b+1,c+1,x)$ decreases in $c$, hence   (\ref{der}) implies that $D_x F(a,b,c,x)$ increases in $c$; together with $F(a,b,c,0)=1$ 
this gives  that $F(a,b,c,x)$ increases in $c$.
Next,  for the generic $b<0$, from a contingency relation    (\cite{Bateman} p. 104 Equation (43))
$$F(a,b-1,c,x)=(1-x)F(a,b,c,x)+(1-a/c)x \,F(a,b,c+1,x),$$
follows that $F(a,b-1,c,x)$ increases in $c$ provided $F(a,b,c,x)$ does.
\end{proof}

\begin{lemma}\label{L1-HH} For $0<x<1$, we have 
 the sign identities:
$$\sgn D_b \left(\frac{F(a_1,b,c,x)}{F(a_2, b,c,x)}\right)=\sgn (a_1-a_2) {\rm ~~if~~}c>b>0,$$
$$\sgn D_a \left(\frac{F(a,b_1,c,x)}{F(a, b_2,c,x)}\right)=\sgn (b_1-b_2)
{\rm ~~if~~}c>a>0$$
and 
$$\sgn D_c \left(\frac{F(a_1,b,c,x)}{F(a_2, b,c,x)}\right)=\sgn (a_2-a_1) {\rm ~~if~~}c>b>0.$$

\begin{proof}  
Applying the integral representation (\ref{Euler}) to 
$$D_b \left(\frac{F(a_1,b,c,x)}{F(a_2,b,c,x)}\right)   =\frac{ F(a_2,b,c,x)   D_b F(a_1,b,c,x)-F(a_1, b,c,x) D_b F(a_2,b,c,x)}{(F(a,b,c,x))^2},$$
the gamma factors cancel and this
has the same sign as

\begin{eqnarray*}
\int_0^1 \frac{z^{b-1}(1-z)^{c-b-1}}{(1-x z)^{a_2}}
{\rm d}z 
\int_0^1 \frac{y^{b-1}(1-y)^{c-b-1}}{(1-x y)^{a_1}} \log\left( \frac{y}{1-y}\right)
{\rm d}y  -      ~~~~~~~~~~~~~~~~~~~~~~~~~~~~~\\
\!\!\!\!\!
\int_0^1 \frac{y^{b-1}(1-y)^{c-b-1}}{(1-x y)^{a_1}} 
{\rm d}y
\int_0^1 \frac{z^{b-1}(1-z)^{c-b-1}}{(1-x z)^{a_2}}\log\left( \frac{z}{1-z}\right)
{\rm d}z 
=~~~~~~~~~\\
\iint\limits_{[0,1]^2}
\left\{
\frac{(yz)^{b-1}((1-y)(1-z))^{c-b-1}}{((1-xy)(1-xz))^{a_1}} \right\}
                \left[  (1-xy)^{a_2-a_1} \log\left( \frac{y-yz}{z-yz}\right) \right]{\rm d}y{\rm d}z=~~~~~~~~~
\end{eqnarray*}
\begin{eqnarray*}
 \iint\limits_{0<z<   y<1} \left\{ \cdots\right\}  \left[  (1-xy)^{a_2-a_1} \log\left( \frac{y-yz}{z-yz}\right) \right]{\rm d}y{\rm d}z   
+   ~~~~~~~~~~~~~~~~~~~~~~~~~~~~~~~~~~~~~~~~~ \\
\iint\limits_{0<y<   z<1} \left\{ \cdots\right\}  \left[  (1-xy)^{a_2-a_1} \log\left( \frac{y-yz}{z-yz}\right) \right]{\rm d}y{\rm d}z  =~~~~~~~~~\\
 \iint\limits_{0<z<   y<1} \left\{ \cdots\right\}  \left[  (1-xy)^{a_2-a_1} \log\left( \frac{y-yz}{z-yz}\right) \right]{\rm d}y{\rm d}z +   ~~~~~~~~~~~~~~~~~~~~~~~~~~~~~~~~~~~~~~~~~
  \\
\iint\limits_{0<z<  y <1} \left\{ \cdots\right\}  \left[  (1-xz)^{a_2-a_1} \log\left( \frac{z-yz}{y-yz}\right) \right]{\rm d}z{\rm d}y  =~~~~~~~~~\\
~~~~~~~~~~~\iint\limits_{0<z<   y<1} \left\{ \cdots\right\}  \left[(1-xy)^{a_2-a_1} - (1-xz)^{a_2-a_1}   \right]        \log\left( \frac{y-yz}{z-yz} \right)   {\rm d}z{\rm d}y\,. ~~~~~~~~~~~~~~~~~
\end{eqnarray*}
But for $y>z$ 
the last integrand has the same sign as $a_1-a_2$. 

The second sign identity follows from symmetry (\ref{symmetry}), and the third is shown by a similar argument.
\end{proof}

\end{lemma}

\begin{lemma}\label{L2-HH} 
For $c>a>0$ and $0<x<1$
$$\sgn D_a \left(  \frac{  F(a,b,c,x)}{F(a,b,c+1,x)}\right) =\sgn b\,.$$ 
\end{lemma}

\begin{proof} Re-write the contingency relation  found
in \cite{Bateman}  (p.103 Equation (35)) as
$$
\frac{F(a,b,c,x)}{F(a,b,c+1,x)} =1-\frac{b}{c}+\frac{b}{c}\, \frac{F(a,b+1,c+1,x)}{F(a,b,c+1,x)}.
$$
By Lemma \ref{L1-HH}, the derivative in $a$ of the right-hand side has the same sign as $b$.
\end{proof}

If the second parameter is a negative integer, $F$ is a polynomial. In this case there is a formula for $D_a F$ due to Fr{\"o}hlich \cite{Froehlich}:
\begin{lemma}\label{L3-HH}    For $m$ a positive integer,
$$
D_a F(a,-m,c,x)=\left(\sum_{j=0}^{m=1} \frac{1}{a+j}\right) F(a,-m,c,x)- \sum_{j=0}^{m-1} \frac{m!}{j!(m-j)} F(a,-j,c,x).
$$
\end{lemma}

\subsection{Baskakov's operator}

The operator associating with sequence $(u_n,~n\geq0)$ the function
$$f(x)=(1-x)^\nu \sum_{n=0}^\infty \frac{(\nu)_n}{n!} \,x^n\, u_n, ~~~x\in[0,1],~\nu>0$$
is one of  Baskakov's operator from the interpolation theory, see Equation (6) in \cite{Baskakov}.

\begin{lemma}\label{Baskakov} If $(u_n)$ is unimodal then $f(x)$ is unimodal too.
\end{lemma}
\begin{proof}
A lengthy calculation yields the derivative
$$f'(x)=(1-x)^{\nu-1}\sum_{n=0}^\infty \frac{(\nu)_n}{n!}  (\nu+n)    (u_{n+1}-u_n) \,x^n\,.$$
By unimodality of $(u_n)$,  the coefficients of the power series $\sum (\cdots)x^n$ have at most one variation of sign, 
 By Descarte's rule of signs \cite{Descartes} $f'$ has at most one positive root, hence $f$ is unimodal on $[0,1]$.

\end{proof}

\subsection{Nevzorov's model for records}\label{Nevzorov}

This section gives some motivation for (\ref{profile}) by way of the theory of records. 
Recall 
the following well known fact.  For i.i.d. random variables $X_1, X_2,\cdots$ with continuous distribution, 
 the indicators of maximal records
$I_k:=1_{\{X_k=\max(X_1,\cdots,X_k)\}}$
are independent and satisfy
${\mathbb P}(I_k=1)=1/k$. 
Nevzorov's model is a generalisation to the nonstationary case.

For  continuous distribution function $G$, let $X_k$  be independent with 
$${\mathbb P}(X_k\leq x)=(G(x))^{\gamma_k},~~~k=1,2,\cdots,$$
where  $\gamma_k$ are given positive constants.
Then the indicators are independent with
$${\mathbb P}(I_k=1)=\frac{\gamma_k}{\gamma_1+\cdots+\gamma_k},$$
see \cite{Nevzorov} (Lecture 25).
A special choice of constants yields the profile (\ref{profile}).
\begin{assertion}\label{NM} The profile of indicators $I_1, I_2,\cdots$ is {\rm (\ref{profile})} for 
$$\gamma_k=\frac{(\theta)_{k-1}}{(k-1)!}\,,$$
and this is the only choice of $\gamma_k$'s up to a constant  positive multiple.
\end{assertion}
\begin{proof} This is shown by induction with the aid of the mean-value recursion
$$\gamma_k= \frac{\theta}{k-1}\sum_{i=1}^{k-1}\gamma_i\,.$$\end{proof}

For  $\theta>1$ integer, $\gamma_k=\frac{(\theta)_{k-1}}{(k-1)!}={\theta+k-2\choose k-1}$ is integer too. In that case 
  the last-success problem has  interpretation as the  best-choice problem with `group interviewing', where  
at stage $k$ the choice can be made from $\gamma_k$ items sampled without replacement from a population of rankable individuals.
For  large number of groups this model is asymptotically simiar  to the case $\gamma_k=k^{\theta-1}$ considered in
\cite{Pfeifer}.

Hofman in \cite{Hoffman} applied Nevzorov's model for records to the best choice problem where arrivals are Poissonian with unknown intensity. However, monotonicity has been stated without proof.

\end{document}